\documentclass[12pt,reqno,oneside]{amsproc}
\usepackage[letterpaper,margin=1in,headheight=13pt,headsep=16pt]{geometry}
\usepackage{amsmath,amssymb,amsthm}
\usepackage{esint}
\usepackage{times}
\usepackage{microtype}
\usepackage{fancyhdr}
\usepackage{needspace}
\usepackage{xcolor}
\usepackage{tikz}
\usepackage{enumitem}
\usetikzlibrary{arrows.meta,calc}
\definecolor{mygray}{gray}{0.1}
\definecolor{cole}{rgb}{0.45,0.0,0.0}

\chardef\showllabel=0
\chardef\forshowkeys=0
\chardef\refcheck=0
\chardef\showcolors=0

\ifnum\showcolors=1
\def\cole{\color{cole}}

\def\colb{\color{black}}
\else  
\def\cole{\color{black}}

\def\colb{\color{black}}
\fi

\ifnum\forshowkeys=1

\usepackage[notcite,color]{showkeys}
\definecolor{refkey}{gray}{0.1}
\fi

\usepackage[
colorlinks=true,
linkcolor=blue,
citecolor=blue,
urlcolor=blue
]{hyperref}

\ifnum\showllabel=1
\def\llabel#1{%
\marginnote{\color{mygray}
\rmfamily\small (#1)}[-0.0cm]\notag
}
\else
\def\llabel#1{\notag}
\fi

\ifnum\refcheck=1
\usepackage{refcheck}
\fi
\newtheorem{theorem}{Theorem}[section]
\newtheorem{lemma}[theorem]{Lemma}

\theoremstyle{definition}

\theoremstyle{remark}

\numberwithin{equation}{section}

\newcommand{\R}{\mathbb R}

\newcommand{\Hs}{\mathcal H}
\newcommand{\diver}{\operatorname{div}}
\newcommand{\diam}{\operatorname{diam}}

\newcommand{\Nu}{N_u}
\newcommand{\T}{\mathcal T}
\newcommand{\TT}{\mathcal T_{\mathrm{term}}}

\allowdisplaybreaks[1]

\hypersetup{pdftitle={Nodal sets on Denjoy-Carleman manifolds},pdfsubject={Weight-scale finite-order estimates and adaptive nodal-set summation}}
\title[Hausedorff measure of nodal sets]{Upper bounds for nodal sets of Laplace eigenfunctions\\
	on Denjoy--Carleman manifolds}
\author[L.~Li]{Linfeng Li}
\address{Department of Mathematics, Southern University of Science and Technology, Shenzhen, Guangdong 518055, People's Republic of China}
\email{lilf3@sustech.edu.cn}
\date{}

\begin{document}
	\def\MM{\mathcal{M}}
	\def\NN{\mathbb{N}}
	\def\linfeng#1{\textbf{\color{red}{LL: #1}}}
	\def\coly{\color{gray}}
\def\gdot{{\color{green}{\hskip-.0truecm\rule[-1mm]{4mm}{4mm}\hskip.2truecm}}\hskip-.3truecm}
\def\comma{ {\rm ,\qquad{}} }  
\def\pa{\partial}
\def\llabel#1{\nonumber}
\begin{abstract}
Let $(M,g)$ be a compact connected Riemannian manifold
without boundary of Denjoy--Carleman regularity.
We establish new upper bounds for the Hausdorff measure of the
nodal sets of Laplace eigenfunctions on $M$.
In particular, for Gevrey manifolds with exponent $s>1$,
we obtain the estimate $C\lambda^{p_s/2}$,
where $\lambda$ is the eigenvalue and $p_s\in (1,s)$ is uniformly bounded for all $s>1$.
We also reduce the logarithmic loss for quasianalytic
classes with logarithmic weights.
These results improve the estimates of
Hezari (\emph{Comm. Math. Phys.}, 2023).
\end{abstract}

\maketitle

\section{Introduction}
\label{sec01}
Let $(M,g)$ be a compact connected Riemannian manifold without
boundary of dimension $n\geq2$, and let $\phi_\lambda$ be a nonzero
real eigenfunction satisfying
\begin{align} 
	-\Delta_g\phi_\lambda=\lambda\phi_\lambda,
		\llabel{EQ00d}
\end{align}
where $-\Delta_g$ is the Laplace--Beltrami operator on $(M,g)$ and
$\lambda>0$ is the eigenvalue.
A fundamental question is how the nodal (zero) set of $\phi_\lambda$ depends on the eigenvalue
and the regularity of the metric.
Yau's conjecture~\cite{Y} predicts that
\begin{align}
	C_1\sqrt\lambda\leq\Hs_g^{n-1}(\{x\in M:\phi_\lambda(x)=0\})
	\leq 
	C_2\sqrt\lambda,
	\llabel{EQ00c}
\end{align}
where $C_1, C_2>0$ are constants depending only on $(M,g)$.
This conjecture was resolved in \cite{DF1, DF3} for real analytic manifold with and without boundary by
Donnelly and Fefferman.
Later, Lin \cite{L} studied nodal sets of elliptic and parabolic equations and obtained sharp upper bound for solutions of the heat equation on real analytic manifolds.
Kukavica \cite{K} extended the results of Donnelly and Fefferman to eigenfunctions of regular elliptic operators of arbitrary even order with real analytic coefficients, establishing sharp upper bounds for their nodal sets.

For smooth manifolds, Hardt and Simon~\cite{HS} proved an exponential-type upper bound for the
nodal sets of Laplace eigenfunctions.
Recently, there were significant progress towards Yau's conjecture for smooth manifolds made by Logunov and Malinnikova in \cite{Lo1, Lo2, LM1}.
For two-dimensional manifolds, Logunov and Malinnikova showed in \cite{LM1} that
\begin{align}
	\Hs_g^{1}(\{x\in M:\phi_\lambda(x)=0\})
	\leq 
	C \lambda^{\frac{3}{4}-\beta},
	\llabel{EQ00a}
\end{align}
for some $\beta \in (0,1/4)$, which improves the upper bound $C\lambda^{3/4}$ by Donnelly and Fefferman \cite{DF2} and Dong \cite{Do}.
In \cite{Lo1}, Logunov proved in higher dimensions $n\geq 3$ that
\begin{align}
	\Hs_g^{n-1}(\{x\in M:\phi_\lambda(x)=0\})
	\leq 
	C \lambda^{\alpha_n},
	\llabel{EQ00b}
\end{align}
where $\alpha_n>1/2$ is a constant depending only on the dimension.
The proof combines the almost monotonicity of the doubling index, quantitative propagation of smallness, and a combinatorial argument that controls the number of subcubes with large doubling indices.
In \cite{Lo2}, Logunov resolved the Nadirashvili's conjecture and established the lower bound in Yau's conjecture for smooth compact Riemannian manifolds without boundary.
We refer to~\cite{LM2} for a survey of progress on Yau's conjecture and to \cite{CZ, GR, KL, LZ, KZZ, LMNN, LS, LTY, NV, Z1, Z2, ZZ} for related results on the nodal sets of elliptic equations in various settings.
Hezari~\cite{He} established upper bounds for nodal sets of Laplacian eigenfunctions on ultradifferentiable Riemannian manifolds, including
\begin{align}
	\Hs_g^{n-1}(\{x\in M:\phi_\lambda(x)=0\})
	\leq C\lambda^{s/2}
	\llabel{EQ00e}
\end{align}
for Gevrey manifolds of order $s\geq1$.
However, this bound deteriorates as the Gevrey exponent increases and eventually becomes weaker than Logunov's polynomial upper bound~\cite{Lo1}.
This motivates us to seek sharper estimates whose exponents remain uniformly bounded as the Gevrey exponent increases.
In related work, Ignatova and Kukavica~\cite{IK1,IK2}
established quantitative unique continuation and Hausdorff measure of nodal sets
for higher-order parabolic and elliptic equations with Gevrey coefficients in one and two spatial dimensions, respectively.
Camliyurt, Kukavica, and Li~\cite{CKL} obtained upper bounds for the spatial nodal sets of parabolic equations with Gevrey-regular coefficients in any dimension, using quantitative unique continuation and interpolation along lines.

The goal of this paper is to improve the upper bounds of nodal sets obtained by Hezari~\cite{He} for ultradifferentiable Riemannian manifolds belonging to Denjoy--Carleman classes, which lie between the real-analytic and smooth classes.
Our estimates depend explicitly on the derivative growth
of the metric and yield sharper bounds on the nodal measure.
In particular, we obtain a smaller exponent of the eigenvalue
in the Gevrey case and reduce the logarithmic loss for
quasianalytic manifolds with logarithmic weights.
These results are stated in Theorems~\ref{T01} and~\ref{TQA},
respectively, and both follow from the general
Denjoy--Carleman estimate in Theorem~\ref{TDC}.
The main difficulty is to combine the derivative estimates
provided by the regularity class with the local growth measured
by the doubling index, since the local nodal estimate is valid
only on sufficiently small spatial scales depending on the
doubling index.
Another difficulty arises from general Denjoy--Carleman
weights, which need not have the factorial growth as in Gevrey
weights.

Our proof begins by establishing a local nodal estimate for
solutions of the divergence-form elliptic equations with
Denjoy--Carleman coefficients.
Using iterated elliptic estimates,
we obtain analytic-type derivative bounds up to a prescribed
finite order on sufficiently small spatial scales
(cf.~Lemma~\ref{P01}).
Combining these bounds with quantitative lower bounds
derived from the doubling index, we control the number
of zeros along suitably chosen lines through interpolation,
drawing also on ideas from~\cite{IK1}.
Following the geometric argument in~\cite{CKL},
we then obtain a local Hausdorff measure estimate for
the nodal set that is linear in an upper bound for
the doubling index (cf.~Lemma~\ref{P02}).
To obtain a global estimate, we construct an adaptive
stopping tree (cf.~Figure~\ref{tree}) based on Logunov's cube-counting
estimate~\cite{Lo1}, subdividing cubes until the local
nodal estimate becomes applicable.
Finally, a suitable weighted summation over the terminal
cubes, combined with the standard harmonic lifting argument,
yields the desired global estimates for the nodal sets
of Laplace eigenfunctions.

The paper is organized as follows.
In Section~\ref{sec:results}, we introduce the assumptions on the Denjoy--Carleman weight sequence and states the main results, Theorem~\ref{T01}--\ref{TDC}.
In Section~\ref{sec02}, we establish finite-order elliptic
estimates and derive a local Hausdorff measure bound for
nodal sets.
We then construct an adaptive stopping tree and obtain
a quantitative nodal estimate through a weighted summation
over the terminal cubes.
In Section~\ref{sec06}, we apply these local estimates
to the harmonic lift of Laplace eigenfunctions and prove
the general Denjoy--Carleman theorem, from which the
Gevrey and quasianalytic results follow.

\section{Main results}
\label{sec:results}
The Denjoy--Carleman classes provide a common framework for
analytic, Gevrey, and more general ultradifferentiable
regularity. 
They prescribe the growth of derivatives
through a weight sequence.
Our aim is to express the nodal set estimates directly
in terms of this sequence, so that different regularity
classes can be treated within the same argument.
We recall the definitions needed below.

We assume that
$\mathcal M=(\MM_k)_{k\geq0}$ is a sequence of positive numbers such that
\begin{enumerate}[label=(\roman*)]
\item $\MM_0=1$;
\item $\MM_k$ is nondecreasing;
\item $\MM_k^2\leq\MM_{k-1}\MM_{k+1}$ for every $k\geq1$;
\item $\displaystyle\sup_{k\geq1}
\left(\frac{\MM_{k+1}}{\MM_k}\right)^{1/k}<\infty$.
\end{enumerate}
Throughout the paper we also assume moderate growth, namely,
\begin{enumerate}[label=(\roman*),start=5]
\item there exists $H\geq1$ such that
\begin{equation}\label{EQ105a}
\MM_{j+k}\leq H^{j+k}\MM_j\MM_k,
\qquad j,k\in\mathbb N_0.
\end{equation}
\end{enumerate}
Although (iv) follows from (v) by taking $j=1$,
we keep it to distinguish the regularity assumptions
from the additional moderate-growth condition.
We refer to \cite[Section~2]{He} for further discussion
of Denjoy--Carleman classes and to~\cite{Fu} for their
geometric framework.

Let $\Omega\subset \mathbb{R}^d$ be open.
A function $f\in C^\infty (\Omega)$ belongs to
$C^{\mathcal M}(\Omega)$ if, for every $W\Subset\Omega$,
there exist constants $K,\rho>0$ such that
\begin{align}
	\|\partial^\alpha f\|_{L^\infty(W)}
	\leq K\rho^{-|\alpha|}|\alpha|!\MM_{|\alpha|}
	\comma
	\alpha\in\mathbb{N}_0^d.
	\label{EQ130a}
\end{align}
The choice $\MM_k=1$ gives the real-analytic class and $\MM_k=(k!)^{s-1}$ gives the Gevrey class~$G^s$ for $s>1$, respectively.
The class $C^{\mathcal M}$ is called quasianalytic if
any function in the class whose derivatives of every
order vanish at one point must vanish identically
on the connected component containing that point.
Under the assumptions (i)--(iii), the Denjoy--Carleman theorem states that $C^{\mathcal{M}}$
is quasianalytic if and only if
\begin{equation}\label{EQ105c}
	\sum_{k=1}^{\infty}\frac{\MM_{k-1}}{k\MM_k}=\infty.
\end{equation} 

Let
\begin{equation}\label{EQ103a}
\omega(t):=\MM_{\lceil t\rceil}^{1/\lceil t\rceil},
\qquad t
\geq 
1,
\end{equation}
where $\lceil \cdot \rceil$ is the ceiling function.
It turns out that $\omega$ is a central quantity in our analysis, as it determines the scale at which the local nodal set estimate applies (see Lemmas~\ref{P01}--\ref{L07} for instance).

We say that a Riemannian manifold $(M,g)$ belongs
to $C^{\mathcal M}$ if it admits an atlas whose
transition maps belong to $C^{\mathcal M}$ and
whose metric coefficients belong to the same class
in these coordinates.
The Riemannian manifolds considered in this paper are assumed to be compact, connected, and boundaryless of
dimension $n\geq2$, and all eigenfunctions are nonzero and real.

Set $d=n+1$.
We fix the constants $B_0\geq2$ and $0<c_0<1$,
depending only on $d$, from \cite[Theorem~5.1]{Lo1},
which is restated in Lemma~\ref{L01} below.
Define
\begin{equation}
	\label{EQ134}
	P_d:=\frac{\log(2B_0-1)}{\log(1+c_0)}
\end{equation}
and
\begin{equation}
	\label{EQ134a}
	\theta_d:=\frac{\log\!\left(\frac12+
		(B_0-\frac12)(1+c_0)^{-1}\right)}{\log B_0}.
\end{equation}
Then $P_d>1$ and $0<\theta_d<1$.

Denote the nodal (zero) set by
\begin{align}
	Z_{\phi_\lambda}= \{x\in M: \phi_\lambda (x)=0\}.
	\llabel{EQ134b}
\end{align}

The following three theorems establish new upper bounds for the nodal sets of Laplace eigenfunctions on Riemannian manifolds of Gevrey, quasianalytic, and general Denjoy--Carleman regularity, respectively.
\cole
\begin{theorem}
	\label{T01}
Let $(M,g)$ belong to the Gevrey class $G^s$, where $s>1$.
There exists constants $p_s\in(1,\min\{s,P_d\})$ and $C>0$ such that
for every eigenfunction $\phi_\lambda$ with $\lambda\geq1$, we have
\begin{equation}\label{eq1}
\Hs_g^{n-1}(Z_{\phi_\lambda})\leq C\lambda^{p_s/2}.
\end{equation}
More precisely, $p_s$ is the unique solution
in $(1,s)$ of
\begin{equation}
	\label{EQ105b}
B_0^{(p_s-1)/(s-1)}
=\frac12+\left(B_0-\frac12\right)(1+c_0)^{-p_s}.
\end{equation}
\end{theorem}
\colb

Our result improves the upper bound $\lambda^{s/2}$ in \cite{He} and provides an explicit characterization of the improved exponent $p_s$ as the unique solution of \eqref{EQ105b}.
In particular, $p_s\leq P_d$, where $P_d$ depends only on the dimension of the manifold. Thus, the exponent remains uniformly bounded as $s\to\infty$, consistent with the polynomial upper bounds established in \cite{Lo1,LM1} for smooth manifolds, which also apply to the Gevrey setting.

\cole
\begin{theorem}
	\label{TQA}
	Let $\beta \in (0,1]$.
	Let $(M,g)$ belong to $C^{\mathcal M}$, where $	\MM_k=(\log(k+e))^{\beta k}$.
	Then there exists a constant $C>0$ such that
	for every eigenfunction $\phi_\lambda$ with
	$\lambda\geq1$, we have
	\begin{equation}
		\label{EQ105e}
		\Hs_g^{n-1}(Z_{\phi_\lambda})
		\leq C\sqrt\lambda\,
		(\log(\lambda+e))^{\beta\theta_d}.
	\end{equation}
\end{theorem}
\colb

The weight sequence is
quasianalytic, as verified in Appendix~\ref{App2}.
The estimate \eqref{EQ105e} improves the logarithmic
loss in Hezari's result~\cite{He} and approaches
the upper bound predicted by Yau's conjecture~\cite{Y}
up to a smaller logarithmic factor.

We now state a general estimate for Denjoy--Carleman
manifolds, from which Theorems~\ref{T01} and~\ref{TQA}
follow.
\cole
\begin{theorem}
	\label{TDC}
Let $(M,g)$ belong to $C^{\mathcal M}$, where $\mathcal M$ satisfies
(i)--(v). 
Let $p\geq1$ and $0\leq\theta\leq1$ satisfy
\begin{equation}
	\label{EQ105f}
	B_0^{-\theta}
	\left(
	\frac12+
	\left(B_0-\frac12\right)(1+c_0)^{-p}
	\right)
	\leq1.
\end{equation}
Then there exists a constant $C>0$ such that for every eigenfunction $\phi_\lambda$ with $\lambda\geq1$, we have
\begin{equation}
	\label{EQ105g}
\Hs_g^{n-1}(Z_{\phi_\lambda})
\leq C\lambda^{p/2}
\sup_{1\leq t\leq\sqrt\lambda}
\frac{\omega(t)^\theta}{t^{p-1}}.
\end{equation}
In particular, we have
\begin{equation}\label{EQ105h}
\Hs_g^{n-1}(Z_{\phi_\lambda})
\leq C\min\left\{
\sqrt\lambda\,\omega(\sqrt\lambda)^{\theta_d},
\lambda^{P_d/2}\right\}.
\end{equation}
\end{theorem}
\colb

The estimate \eqref{EQ105h} follows from \eqref{EQ105g}
by choosing $(p,\theta)=(1,\theta_d)$ and $(P_d,0)$,
which yield the first and second bounds on the
right-hand side, respectively.
The first bound improves the corresponding estimate
of Hezari~\cite{He} by reducing the dependence on
the weight function $\omega$, while the second
ensures that the exponent of $\lambda$ remains
uniformly bounded, independently of the regularity
class, which is consistent with Logunov's polynomial
upper bound~\cite{Lo1}.

\section{The local equation and the doubling index}
\label{sec02}
In this section, we consider a nonzero real solution of the elliptic equation
\begin{equation}
	\diver(A\nabla u)=0
	\label{EQ01}
\end{equation}
in a bounded domain $\Omega\subset \mathbb{R}^d$, where $d\geq 2$.
We assume that the coefficient matrix $A$ is real symmetric and uniformly elliptic, namely, there exists a constant $\Lambda>0$ such that
\begin{equation}
	\label{EQ02}
	\Lambda^{-1}
	|\xi|^2
	\leq A(x)\xi\cdot\xi
	\leq
	\Lambda
	|\xi|^2
	\comma
	\xi \in \mathbb{R}^d
	\comma x\in \Omega.
\end{equation}
Additionally, we assume that $A$ belongs to the Denjoy--Carleman
class associated with the fixed sequence $\mathcal M$ from
Section~\ref{sec:results}.
More precisely, there exist constants $K,\rho>0$ such that
\begin{equation}\label{EQ02b}
\|\partial^\alpha A\|_{L^\infty(\Omega)}
\leq K\rho^{-|\alpha|}|\alpha|!\MM_{|\alpha|},
\qquad \alpha\in\mathbb N_0^d.
\end{equation}
In this section, $u$ is assumed to satisfy
\eqref{EQ01}--\eqref{EQ02b} in $\Omega$, unless stated
otherwise. 
We use only interior cubes in $\Omega$, namely, all balls occurring in the
estimates in this section are required to lie in $B_1=B(0,1)\Subset \Omega$.

For a closed cube $Q$, let $\ell(Q)$ be its side length and $c_Q$ its center.
We define the doubling index
\begin{equation}
\label{EQ03}
\Nu(Q)=\sup_{\substack{x\in Q\\0<r<\diam(Q)}}
\log\frac{\|u\|_{L^\infty(B_{10dr}(x))}}
{\|u\|_{L^\infty(B_r(x))}}.
\end{equation}

\cole
\begin{lemma}[{\cite[Theorem~5.1]{Lo1}}]
\label{L01}
There exist $0<c_0<1$ and an integer $B_0\geq2$ depending only on the
dimension $d$, and positive constants
$R_0=R_0(\Lambda,K,\rho,\MM_1)$ and
$N_0=N_0(\Lambda,K,\rho,\MM_1)\geq1$, such that the following holds.
If a cube $Q\subset B(0,R_0)$ is divided into $B_0^d$ subcubes
$q_i$ of equal size, then
\begin{equation}
	\label{EQ04}
\#\left\{q_i:N_u(q_i)>
\max\left\{\frac{N_u(Q)}{1+c_0},N_0\right\}\right\}
<\frac12B_0^{d-1}.
\end{equation}
\end{lemma}
\colb

This lemma was
established by Logunov for harmonic
functions on smooth Riemannian manifolds. It was subsequently
extended by Georgiev and Roy-Fortin~\cite[Theorem~4.1]{GR}
to solutions of second-order uniformly elliptic
equations with Lipschitz leading coefficients.
A formulation for elliptic equations in a Euclidean domain is also recorded in \cite[Theorem~18]{D}, where the dependence of $R_0$ and $N_0$
on the ellipticity constant and the $L^\infty$ and Lipschitz
bounds of the coefficients is made explicit.
In our setting, these quantities are controlled by
$\Lambda,K,\rho,\MM_1$.

By decreasing $R_0$ if necessary, we may assume that $R_0<1$. 
Decreasing $c_0$ if necessary ensures $c_0<1$
without weakening the conclusion.

We denote $d_{*}:=\left\lfloor\frac{d}{2}\right\rfloor+1$. 
Throughout this section, we use $C\geq 1$ to denote a sufficiently large constant, which may depend on $d, \Lambda, K, \rho, H$ without mention; the value of $C$ may vary from line to line.

The next lemma gives an analytic-type derivative
bound through a prescribed finite order, provided
the spatial scale is sufficiently small.
\cole
\begin{lemma}
\label{P01}
Let $x_0\in B_1$ and $r>0$ satisfy $B_{4r}(x_0)\subset B_1$.
There exists a constant $C_0>0$ such that if
$m\in\mathbb N$ satisfies
\begin{equation}
\label{EQ10}
\omega(m+d_{\ast})\leq\frac{\rho}{r},
\end{equation}
then 
\begin{equation}
\|\partial_e^m u\|_{L^\infty(B_r(x_0))}
\leq 
C_0^{m+1}m!\,r^{-m}
\|u\|_{L^\infty(B_{4r}(x_0))}
\label{EQ11}
\end{equation}
for every unit vector $e\in \mathbb{R}^d$.
We emphasize that the constant $C_0>0$ is independent of $m$, $x_0$, and $r$.
\end{lemma}
\colb

\begin{proof}
Denote 
$A_r(y)=A(x_0+ry)$ and $v(y)=u(x_0+ry)$. 
It is easy to see that
\begin{align}
\diver(A_r\nabla v)=0\quad\text{in }B_4,
\label{EQ11a}
\end{align}
Denote 
\begin{align}
	q:=m+d_{*}.
	\llabel{EQ107a}
\end{align}
For $\alpha\in \mathbb{N}_0^d$ with
$|\alpha|=j\leq q-1$, we apply $\partial^\alpha$ to \eqref{EQ11a}, obtaining
\begin{align}
\label{EQ13a}
\diver(A_r\nabla\partial^\alpha v)=-\diver F_\alpha,
\end{align}
where
\begin{align}
	F_\alpha
	:=
	\sum_{0<\beta\leq\alpha}
	\binom\alpha\beta
	\partial^\beta A_r
	\nabla
	\partial^{\alpha-\beta}v.	\label{EQ13b}
\end{align}
Define the radii
\begin{align}
\llabel{EQ12}
R_j=3-\frac jq
\comma
0\leq j\leq q.
\end{align}
It is easy to see that
\begin{align*}
	\Lambda^{-1} |\xi|^2 
	\leq 
	A_r (y) \xi \cdot \xi \leq \Lambda |\xi|^2
	\comma
	\xi \in\mathbb{R}^d
	\comma 
	y\in B_4.
\end{align*}
Therefore, we apply the Caccioppoli inequality to \eqref{EQ13a}, yielding
\begin{align}
\|\nabla\partial^\alpha v\|_{L^2(B_{R_{j+1}})}
\leq
Cq\|\partial^\alpha v\|_{L^2(B_{R_j})} 
+
C \Vert F_\alpha  \Vert_{L^2(B_{R_j})}
\comma 
0\leq j \leq q-1.
\label{EQ14}
\end{align}

For $0\leq j\leq q$, we denote
\begin{align}
	\llabel{EQ12b}
	T_j: =\sum_{|\alpha|=j}\frac{j!}{\alpha!}
	\|\partial^\alpha v\|_{L^2(B_{R_j})}.
\end{align}
Note that
\begin{align}
T_{j+1}
=\sum_{|\gamma|=j+1}
\frac{(j+1)!}{\gamma!}
\Vert \partial^\gamma v\Vert_{L^2 (B_{R_{j+1}})}
=\sum_{|\alpha|=j}\frac{j!}{\alpha!}
\sum_{i=1}^d
\|\partial^{\alpha+e_i}v\|_{L^2(B_{R_{j+1}})},
\label{EQ12c}
\end{align}
where we used the fact that the coefficient of the term indexed by $|\gamma|=j+1$ is equal to
\begin{align}
\frac{(j+1)!}{\gamma!}=\frac{j!}{\gamma!}\sum_{i=1}^d
\gamma_i=\sum_{i:\gamma_i\geq 1}\frac{j!}{(\gamma-e_i)!}.
	\llabel{EQ37a}
\end{align}
From \eqref{EQ14} and \eqref{EQ12c} it follows that
\begin{align}
\begin{split}
\label{EQ15}
T_{j+1}
\leq
C\sum_{|\alpha|=j}
\frac{j!}{\alpha!}
\|\nabla\partial^\alpha v\|_{L^2(B_{R_{j+1}})}
\leq 
C \sum_{|\alpha|=j}
\frac{j!}{\alpha!}
(q \Vert\partial^\alpha v\Vert_{L^2 (B_{R_{j}})}
+
\Vert F_\alpha \Vert_{L^2 (B_{R_j})}
).
\end{split}
\end{align}
From \eqref{EQ13b} and the H\"older's inequality it follows that
\begin{align*}
\begin{split}
&\sum_{|\alpha|=j}
\frac{j!}{\alpha!}
\|F_\alpha\|_{L^2(B_{R_j})}
\leq 
\sum_{|\alpha|=j}
\sum_{0<\beta\leq \alpha}
\frac{j!}{\alpha!}
\binom{\alpha}{\beta}
\Vert \partial^\beta
A_r
\Vert_{L^\infty (B_{R_j})}
\Vert \nabla \pa^{\alpha-\beta}
v\Vert_{L^2 (B_{R_j})}
\\&\qquad
=
\sum_{|\alpha|=j}
\sum_{k=1}^{j}
\sum_{|\beta|=k, \beta\leq \alpha}
\frac{j!}{\alpha!}
\binom{\alpha}{\beta}
\Vert \partial^\beta
A_r
\Vert_{L^\infty (B_{R_j})}
\Vert \nabla \pa^{\alpha-\beta}
v\Vert_{L^2 (B_{R_j})}.
\end{split}
\end{align*}
For $1\leq|\alpha|=j\leq q$, we infer from~\eqref{EQ02b}
and Lemma~\ref{L02} that
\begin{align}
\|\partial_y^\alpha A_r\|_{L^\infty(B_4)}
&\leq K(r/\rho)^j j!\MM_j
=Kj!\bigl((r/\rho)\mu_j\bigr)^j
\leq Kj!\bigl((r/\rho)\omega(q)\bigr)^j
\leq Kj!,
\llabel{EQ106a}
\end{align}
where $\mu_j=\MM_j^{1/j}$, and the last inequality uses~\eqref{EQ10}.
Therefore, using the identities
\begin{equation}
	\llabel{EQ16a}
	\frac{j!}{\alpha!}\binom\alpha\beta
	=\binom jk\frac{k!}{\beta!}
	\frac{(j-k)!}{(\alpha-\beta)!},
\end{equation}
and
\begin{equation}
	\llabel{EQ16b}
	\sum_{|\beta|=k}
	\frac{k!}{\beta!}=d^k,
\end{equation}
we obtain
\begin{align}
\begin{split}
	\sum_{|\alpha|=j}
	\frac{j!}{\alpha!}
	\Vert F_\alpha\Vert_{L^2 (B_{R_j})}	
	&\leq 
	C
	\sum_{k=1}^j
	\binom{j}{k}
	\sum_{|\gamma|=j-k}
	\frac{(j-k)!k!d^k}{\gamma!}
	\Vert \nabla \pa^\gamma v\Vert_{L^2 (B_{R_j})}
	\\&
	\leq
	C \sum_{k=1}^{j}
	\binom{j}{k}
	d^k k! T_{j-k+1},
	\label{EQ15b}
\end{split}
\end{align}
where we used the same reindexing as in \eqref{EQ12c}, since $R_{j-k+1}\geq R_j$.
Combining \eqref{EQ15} and \eqref{EQ15b}, we arrive at
\begin{equation}
	\label{EQ17}
	T_{j+1}\leq CqT_j
	+C\sum_{k=1}^j
	\binom{j}{k} 
	d^k k!T_{j-k+1},
\end{equation}
for any $0\leq j\leq q-1$.

We choose $H_0>1$ sufficiently large that
\begin{equation}
	\label{EQ18}
	\frac{C}{H_0}
	+
	C\sum_{k\geq1}
	\left(\frac{d}{H_0}
	\right)^k
	\leq
	1,
\end{equation}
where $C>0$ is the constant on the right-hand side of \eqref{EQ17}.
We claim that
\begin{equation}
	\label{EQ19}
	T_j
	\leq
	(H_0q)^jT_0,
\end{equation}
for any $0\leq j\leq q$.
The claim is immediate for $j=0$. Assuming it through order $j$ and
using 
\begin{align*}
	\binom jk k!\leq q^k
	\comma
	k=1,2,\ldots,j,
\end{align*}
we obtain
\begin{align*}
	T_{j+1}
	&\leq Cq(H_0q)^jT_0
	+
	C\sum_{k=1}^j
	d^k q^k (H_0q)^{j-k+1}T_0
	\\&
	\leq
	(H_0q)^{j+1}T_0
	\left(\frac{C}{H_0}
	+
	C\sum_{k=1}^j
	\left(\frac{d}{H_0}\right)^k
	\right)
	\leq(H_0q)^{j+1}T_0,
\end{align*}
where we used \eqref{EQ18} in the last step.
Hence, the claim \eqref{EQ19} is proved by induction.
Consequently, for any $0\leq j\leq q$, we have
	\begin{align}
		\sum_{|\nu|=j}
		\|\partial^\nu v\|_{L^2(B_2)}
		&\le
		T_j 
		\le (H_0q)^jT_0.
			\label{EQ50a}
	\end{align}

Fix $\alpha\in \mathbb{N}^d_0$ with $|\alpha|=m$, we use the Sobolev inequality to obtain
\begin{equation}
	\|\partial^\alpha v\|_{L^\infty(B_2)}
	\leq
	C
	\sum_{|\beta|\le d_{*}}
	\|\partial^{\alpha+\beta}v\|_{L^2(B_2)}
	=
	C \sum_{k=0}^{d_{*}}
	\sum_{|\beta|=k}
	\Vert \pa^{\alpha+\beta} v\Vert_{L^2 (B_2)}
	\leq
	C (H_0q)^q T_0
\llabel{EQ53}
\end{equation}
where we used \eqref{EQ50a} in the last step.
Therefore, we have
\begin{equation}
\max_{|\alpha|=m}
\|\partial^\alpha v\|_{L^\infty(B_2)}
\leq
C(H_0q)^qT_0.
\label{EQ54}
\end{equation}

Let $e=(e_1,\ldots,e_d)\in\mathbb R^d$ be a unit vector.
The multinomial theorem gives
\begin{equation}
\partial_e^m v
=
\sum_{|\alpha|=m}
\frac{m!}{\alpha!}\,
e^\alpha\partial^\alpha v.
\llabel{EQ55}
\end{equation}
Therefore, we appeal to \eqref{EQ54} to yield
\begin{align}
\|\partial_e^m v\|_{L^\infty(B_2)}
&\le
\sum_{|\alpha|=m}
\frac{m!}{\alpha!}
|e^\alpha|
\|\partial^\alpha v\|_{L^\infty(B_2)} \le
C (H_0q)^qT_0
\sum_{|\alpha|=m}
\frac{m!}{\alpha!}|e^\alpha|
\leq
C d^{m/2}(H_0q)^qT_0,
\llabel{EQ56}
\end{align}
where we used the Cauchy-Schwarz inequality in the last step.
Using $q=m+d_{*}$ and the Stirling's formula, we deduce that
\begin{align}
	\Vert \pa^m_e v\Vert_{L^\infty (B_2)}
	\leq
	C^{m+1} m! 
	\Vert v\Vert_{L^2 (B_3)}
	\leq
	C^{m+1} m! 
	\Vert v\Vert_{L^\infty (B_3)}.
	\label{EQ56a}
\end{align}

Finally, recalling that $v(y)= u(x_0+ry)$, we conclude the proof of \eqref{EQ11} from \eqref{EQ56a}.
\end{proof}

The following lemma establishes the linear nodal estimate for an admissible cube.
\cole
\begin{lemma}\label{P02}
There exist constants $\tau_0,C>0$ such that the following holds. Let $Q\subset B(0,R_0)$ be a closed cube. 
If there exists a constant $b>0$ with
\begin{equation}
	\label{EQ22a}
b\geq\max\{1,\Nu(Q)\}
\end{equation}
and
\begin{equation}
	\label{EQ22b} \ell(Q)\omega(b)
	\leq
	\tau_0,
\end{equation}
then every zero of $u$ in $Q$ has finite order of vanishing and
\begin{equation}
\label{EQ21}
\Hs^{d-1}(Z_u\cap Q)\leq Cb\,\ell(Q)^{d-1}.
\end{equation}
\end{lemma}
\colb

\begin{proof}
For notational convenience, we write $\ell=\ell (Q)$, $c=c(Q)$, and $R=\sqrt{d} \, \ell$.
Denote
\begin{align}
	\llabel{EQ37b}
	M:=\|u\|_{L^\infty(B_{4R}(c))}>0.
\end{align} 
We divide $Q$ into $(4C_0)^d$ subcubes $q$ of side $a:=\ell/4C_0$, where we have increased the constant $C_0$ from Lemma~\ref{P01} to an integer at least $2$.
For a subcube $q$ with center $c_q$, set
\begin{align}
	x_j^\pm=c_q\pm\frac a4e_j,\qquad j=1,\ldots,d.
	\llabel{EQ36a}
\end{align}
Let $J\geq 1$ be an integer such that $(10d)^{-J}\leq 1/(128C_0\sqrt d)$.
Let $j\in \{1,2,\ldots,d\}$.
We iterate
\eqref{EQ03} at $x_j^\pm$ over the radii
$\ell/(10d)^J,\ell/(10d)^{J-1},\ldots,\ell$, obtaining
\begin{equation}
	\llabel{EQ22}
	\|u\|_{L^\infty(B_{\ell/(10d)^J}(x_j^\pm))}
	\geq e^{-(J+1)b}M,
\end{equation}
where we used $b\geq N_u (Q)$ and $B_{4R}(c)\subset B_{10d\ell}(x_j^\pm)$.
It follows that there exist points
$p_j^\pm\in B_{\ell/(128C_0\sqrt d)}(x_j^\pm)$
such that
\begin{equation}
	\label{EQ23}
	|u(p_j^\pm)|
	\geq
	\frac12e^{-(J+1)b}M.
\end{equation}
It is clear that the points $p_j^{\pm}$ lie in $q$.

We choose a constant $P\geq 1$ sufficiently large such that
\begin{equation}
\label{EQ24}
P\log4>J+1+\log(2C_0)+1.
\end{equation}
We denote
\begin{equation}
\label{EQ25}
m:=\lceil Pb\rceil
\end{equation}
and
\begin{equation}
	\llabel{EQ25a}
	\tau_0:=\frac{\rho}
	{\sqrt d\,H^{\lceil\log_2(P+1+d_*)\rceil}},
\end{equation}
where $H\geq 1$ is the constant from \eqref{EQ105a}.
Since $b\geq1$, we have $m+d_*\leq(P+1+d_*)b$.
By Lemma~\ref{L02} and~\eqref{EQ22b}, we infer that
\begin{align}
\omega(m+d_*)
&\leq H^{\lceil\log_2(P+1+d_*)\rceil}\omega(b) 
\leq H^{\lceil\log_2(P+1+d_*)\rceil}\frac{\tau_0}{\ell}
=\frac{\rho}{\sqrt d\,\ell}=\frac{\rho}{R}.
\llabel{EQ106b}
\end{align}
Therefore, Lemma~\ref{P01} gives
\begin{equation}
	\label{EQ26}
	\|\partial_e^m u\|_{L^\infty(B_R(c))}
	\leq C_0^{m+1}m!R^{-m}M,
\end{equation}
for every unit vector $e\in \mathbb{R}^d$.

Suppose that there exists a line through one of the points $p=p_j^\pm$, where $j=1,\ldots, d$, which contains
$m$ distinct points $z_1,\ldots,z_m$ of $Z_u\cap q$. 
Let $e\in\mathbb{R}^d$ be a unit vector representing the direction of the line.
By Lemma~\ref{L05}, we have
\begin{align*}
\begin{split}
		|u(p)|
	&\leq\frac{\prod_{k=1}^m|p-z_k|}{m!}
	\|\partial_e^m u\|_{L^\infty(q)}
	\leq
		\frac{(\ell \sqrt{d})^m}{m! (4C_0)^m}
	\frac{C_0^{m+1} m! M}{R^m}
	= 
	C_0 4^{-m}M
	<\frac12e^{-(J+1)b}M,
\end{split}
\end{align*}
where we also used
\eqref{EQ24} and \eqref{EQ26}.
This contradicts \eqref{EQ23}. Therefore,
every line through each $p_j^\pm$ intersects $Z_u\cap q$ in fewer than
$m$ points, for any $j=1,\ldots,d$.

Let $x\in Z_u\cap q$.
For the sake of contradiction, we assume that
$\partial^\alpha u(x)=0$
for every $\alpha\in\mathbb{N}_0^d$ with $|\alpha| < m$.
Taylor's theorem centered at $x$, applied along the segment from
$x$ to any $p=p^{\pm}_j$, together with \eqref{EQ26}, gives
the same contradictory upper bound for $|u(p)|$. 
Hence, every zero in
$q$ has a nonzero derivative of some order $j<m$.

The arguments in \cite[Claims~1 and~2]{M} show that the zero
	set of a nontrivial real-analytic function on a connected
	open subset of $\mathbb{R}^d$ is contained in a countable
	union of $(d-1)$-dimensional $C^1$ graphs.
For completeness, we reproduce the arguments in our setting.
For $\beta\in\mathbb{N}_0^d$ with
$|\beta|<m$, we define
\begin{align}
	E_\beta
	:=
	\left\{
	y\in q:
	\partial^\beta u(y)=0, 
	\nabla\partial^\beta u(y)\neq0
	\right\}.
	\llabel{EQ107a}
\end{align}
We claim that
\begin{equation}
	\label{EQ106c}
	Z_u\cap q
	\subseteq
	\bigcup_{\beta\in\mathbb{N}_0^d, |\beta|<m}
	E_\beta.
\end{equation}
Let $x\in Z_u\cap q$. By the finite-order property
proved above, there exists a smallest integer
$k\in\{1,\ldots,m-1\}$ such that
$\partial^\gamma u(x)\neq0$ for some
$\gamma\in \mathbb{N}_0^d$ with $|\gamma|=k$.
By the minimality of $k$, we have
\begin{align}
	\partial^\alpha u(x)
	=
	0,
	\llabel{EQ107b}
\end{align}
for every $\alpha\in\mathbb{N}_0^d$ with $|\alpha|<k$.
We take $i\in\{1,\ldots,d\}$ such that
$\gamma_i\geq1$, and set $\beta=\gamma-e_i$,
where $e_i$ is the $i$th standard coordinate vector.
Then $|\beta|=k-1<m$. 
It is clear that
\begin{align}
	\partial^\beta u(x)=0
	\llabel{EQ107c}
\end{align}
and
\begin{align}
	\partial_i\partial^\beta u(x)
	=
	\partial^\gamma u(x)
	\neq
	0.
	\llabel{EQ107d}	
\end{align}
Thus $x\in E_\beta$, proving
\eqref{EQ106c}.

Fix $\beta\in\mathbb{N}_0^d$ with $|\beta|<m$.
For every $x\in E_\beta$, the implicit function theorem,
applied to $\partial^\beta u$, provides an open
neighborhood $U_x\subset\Omega$ of $x$ such that
\begin{align}
	\left\{
	y\in U_x:\partial^\beta u(y)=0
	\right\}
	:=
	\Gamma_x,
	\llabel{EQ37c}
\end{align}
where, after a permutation of the coordinates,
$\Gamma_x$ is the graph of a $C^1$ function of
$d-1$ variables. In particular,
\begin{align}
E_\beta\cap U_x\subseteq\Gamma_x.	
\llabel{EQ108a}
\end{align}
Let $\mathcal B=\{B_r(a):a\in\mathbb Q^d,\,
r\in\mathbb Q,\ r>0\}$ be a countable basis for the
Euclidean topology.
For every $y\in E_\beta$, there exists an open
neighborhood $U_y$ on which the zero set of
$\partial^\beta u$ is a $C^1$ graph.
Since $\mathcal B$ is a countable basis, there exists some $V_y\in\mathcal B$
such that
\begin{align}
	y\in V_y\subset U_y.
	\llabel{EQ37e}
\end{align}
The collection $\{V_y:y\in E_\beta\}$ contains
at most countably many distinct sets and covers
$E_\beta$. Enumerate these distinct sets as
$\{V_j\}_{j\in J}$, where $J$ is countable.
By construction, each $V_j$ is contained in
some neighborhood $U_y$ on which the zero set
of $\partial^\beta u$ is a $C^1$ graph.
Consequently, $E_\beta$ is contained in a
countable union of $C^1$ graphs.
There are only finitely many $\beta\in\mathbb{N}_0^d$ with
$|\beta|<m$. Hence \eqref{EQ106c}
shows that $Z_u\cap q$ is contained in a countable union
of $C^1$ graphs. 
Since each $C^1$ graph is countably
$(d-1)$-rectifiable, we conclude that
$Z_u\cap q$ is countably $(d-1)$-rectifiable.

Finally, we adapt the argument in the proof of \cite[Lemma~3.4]{CKL} and use the fact
that every line through each $p_j^\pm$, where $j=1,\ldots,d$,
intersects $Z_u\cap q$ in fewer than $m$ points, yielding
\begin{align}
	\mathcal{H}^{d-1}(Z_u\cap q)
	\leq
	Cm\left(\frac{\ell(Q)}{4C_0}\right)^{d-1},
	\llabel{EQ26b}
\end{align}
for every subcube $q$.
Summing over all subcubes and appealing to \eqref{EQ25}, we obtain the desired estimate \eqref{EQ21}. 
\end{proof}

\begin{figure}[htbp]
	\centering
\begin{tikzpicture}[
	x=1cm,y=1cm,
	font=\sffamily\small,
	parent/.style={draw=black!80,fill=gray!12,rounded corners=2pt,
		line width=.85pt,minimum width=2.3cm,minimum height=1.08cm,align=center},
	terminal/.style={draw=black!75,fill=white,rounded corners=2pt,
		line width=.85pt,minimum width=2.3cm,minimum height=1.08cm,align=center},
	childcube/.style={draw=black!70,fill=gray!4,rounded corners=1.2pt,
		line width=.75pt,minimum width=.78cm,minimum height=.78cm},
	link/.style={-{Stealth[length=2.3mm,width=1.75mm]},line width=.8pt,draw=black!72},
	dots/.style={font=\rmfamily\Large,text=black!70}
	]
	% Top node
	\node[parent] (root) at (0,3.25) {nonterminal\\cube};
	% First generation
	\node[terminal] (t1) at (-5.60,0.55) {terminal\\cube};
	\node[terminal] (t2) at (-2.80,0.55) {terminal\\cube};
	\node[parent] (n1) at (0,0.55) {nonterminal\\cube};
	\node[parent] (n2) at (5.60,0.55) {nonterminal\\cube};
	\node[dots] at (2.80,0.55) {$\cdots$};
	\foreach \to in {t1,t2,n1,n2}{\draw[link] (root.south) -- (\to.north);}
	% Continuing branches -- smaller, unlabeled subcubes
	\node[childcube] (a1) at (-1.22,-2.10) {};
	\node[childcube] (a2) at (0,-2.10) {};
	\node[childcube] (a3) at (1.22,-2.10) {};
	\node[childcube] (b1) at (4.38,-2.10) {};
	\node[childcube] (b2) at (5.60,-2.10) {};
	\node[childcube] (b3) at (6.82,-2.10) {};
	\foreach \to in {a1,a2,a3}{\draw[link] (n1.south) -- (\to.north);}
	\foreach \to in {b1,b2,b3}{\draw[link] (n2.south) -- (\to.north);}
	\node[dots] at (0,-3.05) {$\cdots$};
	\node[dots] at (5.60,-3.05) {$\cdots$};
\end{tikzpicture}
\caption{An illustration of a tree $\mathcal{T}$}
\label{tree}
\end{figure}
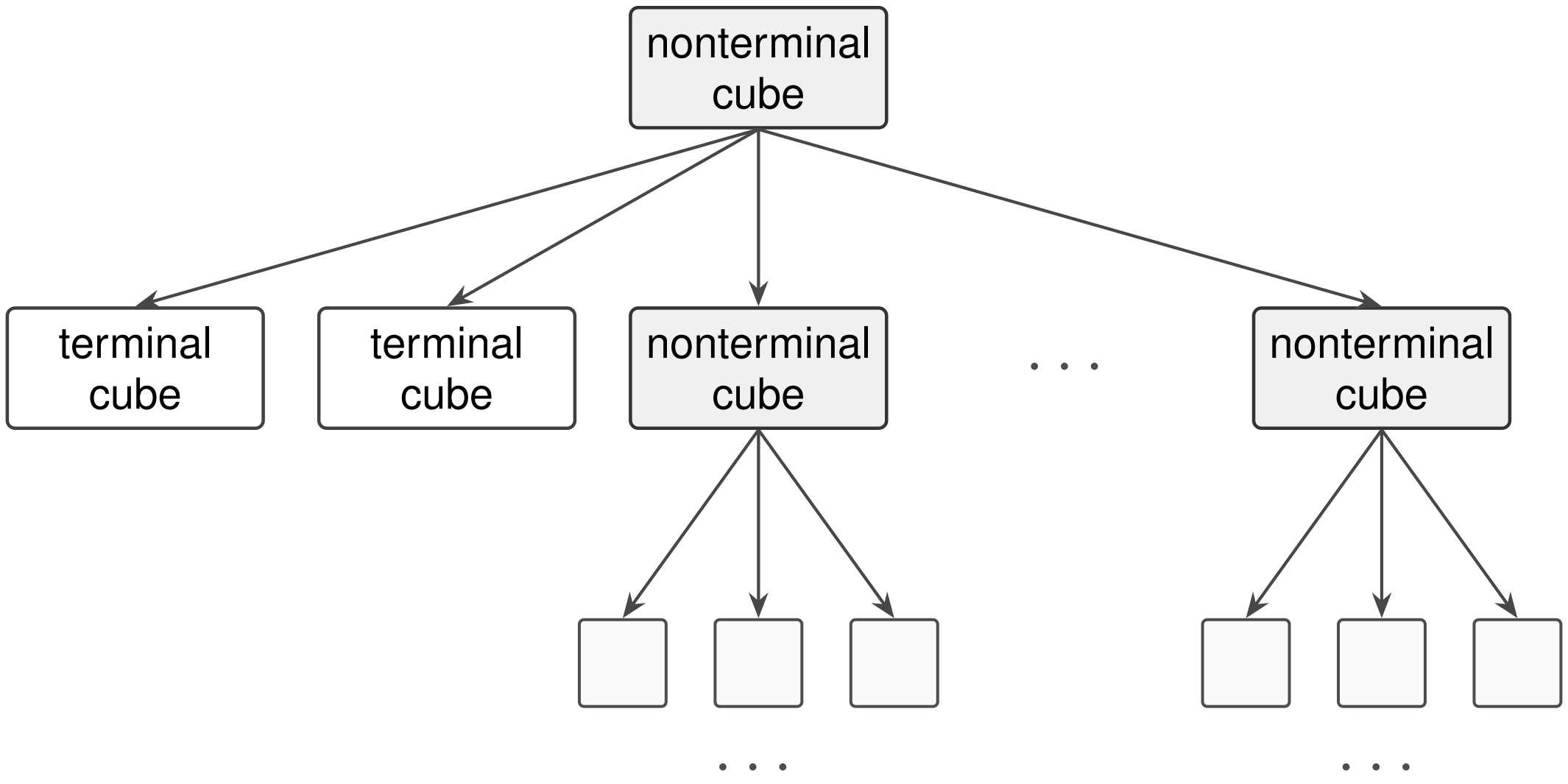

Fix $c_0, B_0, N_0, R_0$ as in Lemma~\ref{L01} and $\tau_0$ as in Lemma~\ref{P02}.
Let $Q_0\subset B(0,R_0)$ be a closed cube and set
\begin{equation}
	\label{EQ32}
	b_0:=\max\{2N_0,\Nu(Q_0)\}.
\end{equation}

We state the setup in the next lemma.
We partition a closed cube into subcubes of equal size. The subdivision proceeds recursively: if a subcube is non-terminal, it is further divided into smaller subcubes, whereas terminal subcubes remain unchanged.
We call the resulting hierarchy of cubes as a tree (see Figure~\ref{tree} for an illustration).
We subdivide until the local nodal estimate in
Lemma~\ref{P02} becomes applicable.
To control this process, we assign to each cube
an upper bound $b(Q)$ for its doubling index.
The assignment permits a fixed-factor decrease
on most children, while preventing an arbitrarily
large decrease in a single generation.
The latter property keeps newly terminal cubes
quantitatively close to the stopping threshold.

We define a function $b$ which assigns a real positive number to a subcube. 
Set
$b(Q_0):=b_0$. 
Given a subcube $Q$ in the tree, we denote
\begin{equation}
\llabel{EQ33}
z(Q)=\frac{\ell(Q)\omega(b(Q))}{\tau_0}.
\end{equation}

In the case when $z(Q)\leq1$, we declare $Q$ terminal. 
This condition is equivalent to
\begin{align}
	\ell(Q)\omega(b(Q))
	\leq
	\tau_0,
	\llabel{EQ33g}
\end{align}
which is precisely
the scale condition \eqref{EQ22b} in Lemma~\ref{P02}.
Together with $b(Q)\geq\max\{1,\Nu(Q)\}$,
which shall be established in the next lemma,
it ensures that Lemma~\ref{P02} is applicable
to every terminal cube.

Otherwise, we divide $Q$ into $B_0^d$
children $q$ of equal size and set
\begin{equation}
	\label{EQ34}
	b(q)=
	\begin{cases}
	b(Q),&\text{~if~}\Nu(q)>b(Q)/(1+c_0),\\
	b(Q)/(1+c_0),&\text{~if~}\Nu(q)\leq b(Q)/(1+c_0).
	\end{cases}
\end{equation}
A child in the first case, i.e., $N_u (q)> b(Q)/(1+c_0)$, is called exceptional. 
We denote the resulting
tree by~$\T$ and the set of terminal subcubes by~$\TT$.

The following lemma establishes some key properties of the tree $\T$, including the control of the doubling index and quantitative bounds on exceptional children.
\cole
\begin{lemma}
\label{L07}
Let $Q_0\subset B(0,R_0)$ be a closed cube with
$\ell(Q_0)\leq\ell_*$, where
\begin{equation}\label{EQ31}
\ell_*:=\min\left\{\frac{R_0}{2},
\frac{\tau_0}{\omega((1+c_0)N_0)}\right\}.
\end{equation}
Then the following statements hold:
\begin{enumerate}
\item Every cube in the tree $\T$ satisfies
\begin{equation}\label{EQ34a}
b(Q)\geq N_u(Q).
\end{equation}
\item Every cube in the tree $\T$ satisfies
\begin{equation}\label{EQ34a1}
b(Q)\geq N_0.
\end{equation}
\item The tree $\T$ has finite depth.
\item Every nonterminal cube has fewer than
$B_0^{d-1}/2$ exceptional children. If $q$ is terminal and its parent
is nonterminal, then
\begin{equation}\label{EQ35}
\frac1{B_0H}<z(q)\leq1,
\end{equation}
where $H\geq 1$ is the constant on the right-hand side of~\eqref{EQ105a}.
\end{enumerate}
\end{lemma}
\colb

\begin{proof}
(1) By~\eqref{EQ32}, we have
\begin{equation}
	\llabel{EQ34b}
b(Q_0)\geq N_u(Q_0).
\end{equation}
Suppose this inequality holds for a nonterminal cube $Q$, and let $q$
be a child. The definition~\eqref{EQ03} gives $N_u(q)\leq N_u(Q)$. 
If $q$ is exceptional, then
\begin{align}
	b(q)=b(Q)\geq N_u(Q)\geq N_u(q).
	\llabel{EQ36b}
\end{align}
Otherwise, the definition~\eqref{EQ34} gives
$b(q)=b(Q)/(1+c_0)\geq N_u(q)$.
This proves~\eqref{EQ34a} by induction.

(2) The assertion is immediate for $Q_0$. Let $q$ have parent $Q$;
then $Q$ is nonterminal. If $b(Q)\leq(1+c_0)N_0$, monotonicity of
$\omega$ from Lemma~\ref{L02} and~\eqref{EQ31} imply
\begin{align}
z(Q)\leq
\frac{\ell_*\omega((1+c_0)N_0)}{\tau_0}\leq1,
	\llabel{EQ109c}
\end{align}
a contradiction. Hence $b(Q)>(1+c_0)N_0$ and
\begin{align}
b(q)\geq\frac{b(Q)}{1+c_0}>N_0.
	\llabel{EQ109d}
\end{align}
This proves~\eqref{EQ34a1}.

(3) For a cube $Q$ at depth $j$, where $j\in\mathbb{N}_0$, we have
$\ell(Q)=\ell(Q_0)B_0^{-j}$ and $b(Q)\leq b_0$. Thus
\begin{align}
z(Q)\leq\frac{\ell(Q_0)\omega(b_0)}{\tau_0B_0^j}
\leq
1,
\llabel{EQ110a}
\end{align}
for $j$ sufficiently large. The tree therefore has finite depth.

(4) If $q$ is an exceptional child of a nonterminal cube $Q$, then
\begin{align}
N_u(q)>\frac{b(Q)}{1+c_0}
\geq\max\left\{\frac{N_u(Q)}{1+c_0},N_0\right\},
	\llabel{EQ110b}
\end{align}
where the last inequality uses (1) and
$b(Q)>(1+c_0)N_0$, which is proven in (2). 
Therefore, Lemma~\ref{L01}  implies that $Q$ has
fewer than $B_0^{d-1}/2$ exceptional children.

For any child $q$ of $Q$, we have
\begin{equation}
	\llabel{EQ106d}
z(q)=\frac{z(Q)}{B_0}
\frac{\omega(b(q))}{\omega(b(Q))}.
\end{equation}
Since $b(Q)/(1+c_0)>N_0\geq1$ and $1<1+c_0<2$,
Lemma~\ref{L02} implies
\begin{align}
\omega(b(Q))
\leq H\omega\!\left(\frac{b(Q)}{1+c_0}\right)
\leq H\omega(b(q)).
	\llabel{EQ110c}
\end{align}
Consequently, $z(q)\geq z(Q)/(B_0H)$. If $q$ is terminal,
$z(Q)>1$ and $z(q)\leq1$, which proves~\eqref{EQ35}.
\end{proof}

The following lemma provides a quantitative upper bound for the
Hausdorff measure of the nodal set and plays a key role
in the proof of Theorem~\ref{TDC}.
\cole
\begin{lemma}
\label{T02}
Let $Q_0\subset B(0,R_0)$ be a closed cube with
$\ell(Q_0)\leq\ell_*$, where $\ell_*$ is given by~\eqref{EQ31}.
If $p\geq1$ and $0\leq\theta\leq1$ satisfy~\eqref{EQ105f},
then
\begin{equation}\label{EQ37}
\Hs^{d-1}(Z_u\cap Q_0)
\leq C\ell(Q_0)^{d-1}b_0^p
\sup_{1\leq t\leq b_0}\frac{\omega(t)^\theta}{t^{p-1}},
\end{equation}
where $C>0$ is a constant independent of $b_0$.
\end{lemma}
\colb

\begin{proof}
In the case when $Q_0$ is terminal, we infer that $z(Q_0)\leq 1$, which leads to
\begin{align}
	\ell (Q_0) \omega (b(Q_0))
	\leq \tau_0.
	\llabel{EQ109a}
\end{align}
From Lemma~\ref{L07} it follows that
\begin{align}
	b_0 = b(Q_0)\geq 
	\max\{N_u (Q_0), N_0 \}
	\geq 
	\max\{N_u (Q_0), 1\}
	\llabel{EQ37a}.
\end{align}
Therefore, the hypotheses of Lemma~\ref{P02} are satisfied and we obtain
\begin{equation}
	\llabel{EQ35a}
\Hs^{d-1}(Z_u\cap Q_0)
\leq Cb_0\ell(Q_0)^{d-1},
\end{equation}
This implies~\eqref{EQ37}, since $\omega\geq1$ and
\begin{align}
	b_0^p\sup_{1\leq t\leq b_0}
	\frac{\omega(t)^\theta}{t^{p-1}}
	\geq b_0\omega(b_0)^\theta\geq b_0.
	\llabel{EQ109b}
\end{align}

Suppose that $Q_0$ is nonterminal. We choose a weight whose total does not increase
under subdivision and which controls the local
nodal contribution on terminal cubes.
More precisely, we define
\begin{equation}
	\label{EQ38}
W(Q):=\tau_0^{-\theta}\ell(Q)^{d-1+\theta}b(Q)^p.
\end{equation}
Then 
\begin{align}
	\ell (Q)^{d-1}
	b(Q) = z(Q)^{-\theta}
	W(Q)
	\frac{\omega (b(Q))^\theta}{b(Q)^{p-1}}.
	\label{EQ38a}
\end{align}
Let $h_Q$ be the number of exceptional children of a nonterminal
cube $Q$. By Lemma~\ref{L07}, we deduce that $h_Q<B_0^{d-1}/2$. 
Therefore, we have
\begin{align}
\begin{split}
	\frac{\sum_{q\text{ child of }Q}W(q)}{W(Q)}
&=B_0^{-(d-1+\theta)}
\left[h_Q+(B_0^d-h_Q)(1+c_0)^{-p}\right]\\
&<B_0^{-\theta}
\left(\frac12+
\left(
B_0-
\frac12
\right)(1+c_0)^{-p}\right)
\leq
1,
\llabel{EQ39}
\end{split}
\end{align}
where we used \eqref{EQ105f} in the last step. 
Since the tree $\T$ has finite depth, successive
summation over all generations of $Q_0$ gives
\begin{equation}\label{EQ40}
\sum_{q\in\TT}W(q)\leq W(Q_0).
\end{equation}

For a subcube $q\in \TT$, from Lemma~\ref{L07} it follows that
\begin{align}
	b(q)
	\geq
	\max\{1,N_u(q)\}.
	\llabel{EQ40e}
\end{align} 
It is clear that $z(q)\leq 1$. Therefore, we infer that
\begin{align}
	\ell(q)\omega(b(q))
	\leq
	\tau_0.
	\llabel{EQ40g}
\end{align}
The hypotheses of Lemma~\ref{P02} are satisfied and it follows that
\begin{align}
\begin{split}
	\Hs^{d-1}(Z_u\cap q)
&\leq Cb(q)\ell(q)^{d-1} =Cz(q)^{-\theta}W(q)
\frac{\omega(b(q))^\theta}{b(q)^{p-1}} 
\\&
\leq 
C(B_0H)^\theta W(q)
\sup_{1\leq t\leq b_0}\frac{\omega(t)^\theta}{t^{p-1}},
\label{EQ40c}
\end{split}
\end{align}
where we used~\eqref{EQ35}, \eqref{EQ38a}, and the fact that $1\leq b(q)\leq b_0$.
The supremum is finite because $\omega$ is a finite-valued
nondecreasing step function on every bounded interval by Lemma~\ref{L02}.

It is easy to see that the terminal subcubes cover $Q_0$. Summing~\eqref{EQ40c} and using
\eqref{EQ38}--\eqref{EQ40}, we obtain
\begin{align}
\begin{split}
	\Hs^{d-1}(Z_u\cap Q_0)
&\leq C(B_0H)^\theta\tau_0^{-\theta}
\ell(Q_0)^{d-1+\theta}b_0^p
\sup_{1\leq t\leq b_0}\frac{\omega(t)^\theta}{t^{p-1}} \\
&\leq C\ell(Q_0)^{d-1}b_0^p
\sup_{1\leq t\leq b_0}\frac{\omega(t)^\theta}{t^{p-1}},
\llabel{EQ40c1}
\end{split}
\end{align}
since $\ell(Q_0)\leq1$. This
proves~\eqref{EQ37} in the case when $Q_0$ is nonterminal.
\end{proof}

\section{Proof of the main results}
\label{sec06}
In this section, we return to the Denjoy-Carleman Riemannian manifold
$(M,g)$ and set $d=n+1$, where $n\geq 2$ is the dimension of the manifold.
We aim to apply the local estimates of Section~\ref{sec02}
to a harmonic lift of the eigenfunction, and then
transfer the resulting estimate back to $M$
using the product structure of the lifted nodal set.
This process is rather standard by now; for completeness, we present the detail below.

We choose finitely many coordinate charts of class $C^{\mathcal{M}}$,
\begin{align}
	\kappa_\nu:\mathcal U_\nu\to\Omega_\nu\subset\R^n
	\llabel{EQ130b}
\end{align}
and open sets $V_\nu\Subset\Omega_\nu$ such that
the sets $\kappa_\nu^{-1}(V_\nu)$ cover $M$.
For each $\nu$, we choose an open neighborhood $W_\nu$ 
such that $V_\nu \Subset W_{\nu} \Subset \Omega_\nu$.
These larger neighborhoods provide room for the enlarged
balls used in the local estimates of Section~\ref{sec02}.
By the positive definiteness of $g$ and the compactness
of $\overline W_\nu$, the eigenvalues of its coordinate
matrix have positive lower bounds and finite upper bounds which are uniform over the finite family of charts.
Moreover, the assumptions (i)--(iii) in
Section~\ref{sec:results} imply that the
Denjoy--Carleman class $C^{\mathcal M}$ is closed under
addition, multiplication, and composition.
Since $\det g$ is bounded away from zero on
$\overline W_\nu$, it follows that $g^{-1}$ and
$w=\sqrt{\det g}$ satisfy bounds of the same form
as~\eqref{EQ130a} on $W_\nu$.
After adjusting the constants, we may choose constants $K,\rho>0$
uniformly over this finite family, without changing
the weight sequence $\mathcal M$.

We fix one chart and suppress the index $\nu$ for convenience.
In these coordinates,
we write $\phi_\lambda(x)$
for $\phi_\lambda(\kappa^{-1}(x))$ and define
\begin{align}
	\llabel{EQ42}
	U_\lambda(x,t)
	=
	e^{\sqrt\lambda\,t}\phi_\lambda(x),
\end{align}
where $(x,t)\in W\times \mathbb{R}.$
Since $V\Subset W$, we may cover $\overline V$
by finitely many sufficiently small closed cubes
$Q\subset W$ such that all enlarged balls used
in the local estimates on $Q\times I$ remain
contained in $W\times\R$, where $I=[-\ell(Q)/2, \ell(Q)/2]$.
In particular, $Q\times I$ is a cube in
$\R^{n+1}$ of side length $\ell(Q)$.

\cole
\begin{lemma}\label{L09}
Let $Q\subset W$ be a sufficiently small cube. Then
\begin{equation}\label{EQ44}
N_{U_\lambda}(Q\times I)\leq C(1+\sqrt\lambda),
\end{equation}
where $C>0$ is a constant independent
of $\lambda$.
\end{lemma}
\colb

The proof of Lemma~\ref{L09} is standard and relies
on the doubling estimate for Laplace
eigenfunctions due to Donnelly and
Fefferman~\cite{DF1}.
 For completeness, we provide the details in Appendix~\ref{App1}.

Now we are ready to prove the main results of the paper.
\begin{proof}[Proof of Theorem~\ref{TDC}]
Before applying the results of Section~\ref{sec02},
we translate and, if necessary, rescale the lifted
coordinate neighborhoods to the normalized setting
considered there.
Since the coordinate atlas is finite, the resulting
ellipticity and Denjoy--Carleman bounds remain uniform.
We may therefore apply Lemma~\ref{L01} and choose
common constants in the local estimates.

We then choose finitely many sufficiently small
coordinate cubes $Q_\nu$ whose inverse images cover $M$.
For each $\nu$, let
\begin{align}
	\llabel{EQ33h}
	I_\nu=
	\left[-\frac{\ell(Q_\nu)}{2},
	\frac{\ell(Q_\nu)}{2}\right].
\end{align}
The cubes are chosen so that, after normalization,
each $Q_\nu\times I_\nu$ satisfies the containment
and side-length conditions of Lemma~\ref{T02}.
These choices are made independently of the eigenvalue~$\lambda$.

Let
\begin{align}
	b_{0,\nu}:=\max\{2N_0,N_{U_\lambda}(Q_\nu\times I_\nu)\}.
	\llabel{EQ114a}
\end{align}
By Lemma~\ref{L09}, we infer that 
\begin{align}
	b_{0,\nu}
	\leq 
	C_*\sqrt\lambda,
	\label{EQ114b}
\end{align}
where $C_*\geq1$ is a constant.

For any constants $T, A, p\geq1$ and
$0\leq \theta \leq 1$, we claim that
\begin{equation}
	\label{EQ106e}
\sup_{1\leq t\leq AT}\frac{\omega(t)^\theta}{t^{p-1}}
\leq H^{\theta\lceil\log_2 A\rceil}
\sup_{1\leq t\leq T}\frac{\omega(t)^\theta}{t^{p-1}}.
\end{equation}
For $t\in [1, T]$ this is immediate.
For $t\in [T, AT]$, we use
$t^{1-p}\leq T^{1-p}$ and
$\omega(t)\leq\omega(AT)\leq H^{\lceil\log_2 A\rceil}\omega(T)$, which follows from Lemma~\ref{L02}.
Hence, \eqref{EQ106e} is proven.

We apply Lemma~\ref{T02} with dimension $d=n+1$ and then
\eqref{EQ106e} with $A=C_*$ and $T=\sqrt\lambda$,
yielding
\begin{equation}
	\label{EQ47}
	\Hs^n
	\bigl(Z_{U_\lambda}\cap(Q_\nu\times I_\nu)
	\bigr)
	\leq C\ell
	(Q_\nu)^n
	\lambda^{p/2}
	\sup_{1\leq t\leq\sqrt\lambda}
	\frac{\omega(t)^\theta}{t^{p-1}},
\end{equation}
where we also appealed to \eqref{EQ114b}.
Since the exponential factor never vanishes, we deduce that
\begin{equation}
	\label{EQ48}
	Z_{U_\lambda}\cap(Q_\nu\times I_\nu)
	=(Z_{\phi_\lambda}\cap Q_\nu)\times 
	I_\nu.
\end{equation}
Since $\ell (Q_\nu) \leq \ell_{*}$, we subdivide the product cube $Q_\nu \times I_\nu$ as in Lemma~\ref{L07} and obtain a finite-depth tree $\T$. 
For every subcube $q\in \TT$, we infer that
\begin{align}
	z(q)=
	\frac{\ell (q) \omega (b (q))}{\tau_0}
	\leq 1.
	\llabel{EQ48a}
\end{align}
Moreover, from Lemma~\ref{L07} it follows that
\begin{align}
	b(q) 
	\geq 
	\max\{N_{U_\lambda} (q), N_0\}
	\geq 
	\max\{N_{U_\lambda} (q), 1\}.
	\llabel{EQ48b}
\end{align}
The hypotheses of Lemma~\ref{P02} are satisfied and we infer that every zero of $U_\lambda$ in $Q_\nu \times I_\nu$ has finite vanishing order.
For the sake of contradiction, we assume that $\phi_\lambda$ has infinite vanishing order at some $x\in Q_\nu$, then
\begin{align}
	\partial_x^\alpha
	\partial_t^jU_\lambda(x,t)
	=
	e^{\sqrt\lambda t}
	\lambda^{j/2}
	\partial_x^\alpha\phi_\lambda(x)=0
	\comma
	\alpha\in\mathbb N_0^n
	\comma 
	j\in\mathbb{N}_0,
	\llabel{EQ48d}
\end{align}
which is a contradiction.
Therefore, every zero of $\phi_\lambda$ has finite vanishing order and
the same stratification argument as in Lemma~\ref{P02}
shows that $Z_{\phi_\lambda}\cap Q_\nu$ is countably
$(n-1)$-rectifiable. 
The product area formula then yields
\begin{equation}
	\label{EQ49}
	\Hs^n
	\bigl((Z_{\phi_\lambda}\cap Q_\nu)
	\times I_\nu\bigr)
	=
	|I_\nu|
	\Hs^{n-1}
	(Z_{\phi_\lambda}
	\cap Q_\nu).
\end{equation}
Combining~\eqref{EQ47}--\eqref{EQ49}, we obtain
\begin{equation}
	\label{EQ50}
	\Hs^{n-1}(Z_{\phi_\lambda}\cap Q_\nu)
	\leq C\ell(Q_\nu)^{n-1}\lambda^{p/2}
	\sup_{1\leq t\leq\sqrt\lambda}
	\frac{\omega(t)^\theta}{t^{p-1}}.
\end{equation}

Finally, on each compact coordinate region, the Riemannian and
Euclidean metrics are uniformly equivalent. Therefore, we have
\begin{align}
	C^{-1}
	\Hs^{n-1}(E)
	\leq\Hs_g^{n-1}(E)
	\leq 
	C\Hs^{n-1}(E),
	\llabel{EQ50a}
\end{align}
for every countably $(n-1)$-rectifiable set $E$ in that region.
Summing~\eqref{EQ50} over the finite cover of the manifold proves~\eqref{EQ105g}.

By~\eqref{EQ134a}, condition~\eqref{EQ105f}
holds with equality for $(p,\theta)=(1,\theta_d)$.
The monotonicity of $\omega$ from Lemma~\ref{L02} then gives the first bound in
\eqref{EQ105h}. For $(p,\theta)=(P_d,0)$, the condition \eqref{EQ105f} again
holds with equality by \eqref{EQ134}.
Since
\begin{align}
	\sup_{1\leq t\leq\sqrt\lambda}t^{1-P_d}=1,
	\llabel{EQ50b}
\end{align}
we obtain the second bound on the right-hand side of \eqref{EQ105h} and thus prove~\eqref{EQ105h}.
\end{proof}

\begin{proof}[Proof of Theorem~\ref{T01}]
Take $\MM_k=(k!)^{s-1}$. It is easy to see that this sequence satisfies the assumptions (i)--(iv) in Section~\ref{sec:results} by direct verification.
The moderate growth assumption (v) follows from
\begin{align}
	(j+k)!
	\leq
	2^{j+k}j!k!
	\comma
	j,k\in\mathbb{N}_0.
	\llabel{EQ51a}
\end{align}

Define
\begin{align}
	F(p):=B_0^{(p-1)/(s-1)}-\frac12-
	\left(B_0-\frac12\right)(1+c_0)^{-p},
	\llabel{EQ115a}
\end{align}
where $B_0\geq 2$ and  $c_0\in (0,1)$ are the constants from Lemma~\ref{L01}.
It is easy to see that the function $F: [1,\infty)\to \mathbb{R}$ is continuous and strictly increasing.
Direct calculation gives $F(1)<0$ and $F(s)>0$.
Therefore, there exists a unique $p_s\in(1,s)$ such that $F(p_s)=0$.
Set
\begin{equation}
	\llabel{EQ29}
\theta=\frac{p_s-1}{s-1}\in(0,1).
\end{equation}
Then
\begin{equation}
	\label{EQ30}
	B_0^{-\theta}
	\left(
	\frac12
	+
	\left(B_0-\frac12\right)(1+c_0)^{-p_s}
	\right)=1.
\end{equation}
For any $t\geq 1$, we have
\begin{equation}
	\llabel{EQ106f}
	\omega(t)
	=
	(\lceil t\rceil!)^{(s-1)/\lceil t\rceil}
	\leq 
	\lceil t\rceil^{s-1}\leq(2t)^{s-1},
\end{equation}
which leads to
\begin{align}
\sup_{1\leq t\leq\sqrt\lambda}
\frac{\omega(t)^\theta}{t^{p_s-1}}
\leq2^{(s-1)\theta}=2^{p_s-1}.
	\llabel{EQ115c}
\end{align}
Therefore, we apply Theorem~\ref{TDC} with $p=p_s$ to obtain~\eqref{eq1}.

Finally, we note that \eqref{EQ30} implies
\begin{align}
\left(B_0-\frac12\right)(1+c_0)^{-p_s}
=B_0^\theta-\frac12>\frac12,
	\llabel{EQ115d}
\end{align}
from where
\begin{align}
p_s<\frac{\log(2B_0-1)}{\log(1+c_0)}=P_d.
	\llabel{EQ115e}
\end{align}
This proves that $p_s\in (1,\min \{s, P_d\})$.
\end{proof}

\begin{proof}[Proof of Theorem~\ref{TQA}]
The logarithmic weights $\mathcal{M}_k= (\log (k+e))^{\beta k}$, where $0<\beta \leq 1$, satisfies the assumptions (i)--(v) and \eqref{EQ105c} by Appendix~\ref{App2}.
Since
\begin{align}
\omega(t)=(\log(\lceil t\rceil+e))^\beta
\leq C_\beta(\log(t+e))^\beta,
\qquad t\geq1,
	\llabel{EQ111d}
\end{align}
we obtain \eqref{EQ105e} by substituting $t$ with $\sqrt\lambda$ in \eqref{EQ105h}.
\end{proof}

\appendix
\section{Auxiliary lemmas}
\label{App}
\cole
\begin{lemma} 
	\label{L05}
	Let $f\in C^m(I)$, where $I\subset\R$ is a compact interval and $m\in \mathbb{N}$.
	If $t_1,\ldots,t_m\in I$ are distinct zeros of $f$, then for every
	$t\in I$, we have
	\begin{equation}
		\label{EQ20}
		|f(t)|\leq\frac{\prod_{k=1}^m|t-t_k|}{m!}
		\|f^{(m)}\|_{L^\infty(I)}.
	\end{equation}
\end{lemma}
\colb
The proof uses the Rolle's theorem which can be found in \cite[Lemma~4.a]{Hi}.
This inequality was also proved via the Hermite intepolation theorem in \cite[Theorem~2.4]{IK1}.

\cole
\begin{lemma}\label{L02}
We assume that the conditions (i)--(v) in Section~\ref{sec:results} hold.
Then the sequence $\{\MM_k^{1/k}\}_{k\geq1}$ and the function
$\omega$ in~\eqref{EQ103a} are nondecreasing.
Moreover, for every $A\geq1$, we have
\begin{equation}\label{EQ101b}
\omega(At)\leq H^{\lceil\log_2 A\rceil}\omega(t),
\qquad t\geq1,
\end{equation}
where $H\geq 1$ is the constant on the right-hand side of~\eqref{EQ105a}.
\end{lemma}
\colb

\begin{proof}
For $k\geq1$, let $\mu_k=\MM_k^{1/k}$ and
$\eta_k=\MM_k/\MM_{k-1}$. Logarithmic convexity implies
$\eta_k\leq\eta_{k+1}$. Since $\MM_0=1$, we obtain
\begin{align}
\mu_k=\left(\prod_{j=1}^k\eta_j\right)^{1/k}
\leq\eta_k\leq\eta_{k+1}.
	\llabel{EQ112b}
\end{align}
It follows that
$\MM_{k+1}=\MM_k\eta_{k+1}\geq\MM_k^{1+1/k}$ and hence
$\mu_{k+1}\geq\mu_k$. The asserted monotonicity of~$\omega$ follows
from~\eqref{EQ103a}.

Moderate growth gives
\begin{align}
\MM_{2k}\leq H^{2k}\MM_k^2.
	\llabel{EQ112c}
\end{align}
Taking the $2k$-th root yields
\begin{equation}\label{EQ101a}
\mu_{2k}\leq H\mu_k,\qquad k\geq1.
\end{equation}
Fix $t\geq 1$.
Let $r=\lceil\log_2 A\rceil$ and $k=\lceil t\rceil$.
Since $At\leq2^rk$ and $2^rk$ is an integer, we obtain 
$\lceil At\rceil\leq2^rk$. 
By the monotonicity of $\omega$ and $r$ iterations
of~\eqref{EQ101a}, we arrive at
\begin{align}
\omega(At)=\mu_{\lceil At\rceil}
\leq\mu_{2^rk}\leq H^r\mu_k=H^r\omega(t).
	\llabel{EQ112d}
\end{align}
For $r=0$, this inequality obviously holds.
Therefore, we conclude the proof of \eqref{EQ101b}.
\end{proof}

\subsection{Proof of Lemma~\ref{L09}}
\label{App1}
	Let
	\begin{equation}\llabel{EQ42a}
		A(x)=w(x)g^{-1}(x)
	\end{equation}
	and
	\begin{equation}
		\llabel{EQ42b}
		\widetilde A(x,t)=
		\begin{pmatrix}A(x)&0\\0&w(x)\end{pmatrix},
	\end{equation}
where $w(x)=\sqrt{\det g}$ and $(x,t)\in W\times \mathbb{R}$.
	It is readily checked that
	\begin{equation}
		\llabel{EQ45}
		-\diver_x(A\nabla\phi_\lambda)
		=
		\lambda w\phi_\lambda,
	\end{equation}
which leads to
	\begin{align*}
		\diver_{x,t}(\widetilde A\nabla_{x,t}U_\lambda)
		&=e^{\sqrt\lambda t}\diver_x(A\nabla\phi_\lambda)
		+w\partial_t^2U_\lambda =-\lambda wU_\lambda+\lambda wU_\lambda=0.
	\end{align*}
	The matrix $\widetilde A$ is real symmetric and uniformly elliptic.
	Moreover, the Denjoy-Carleman bounds for $A$ and $w$ imply that
	 there exist constants $K,\rho>0$ such that
	\begin{align}
	\|\partial_{x,t}^{\alpha}\widetilde A\|_{L^\infty(W\times \mathbb{R})}
	\leq K\rho^{-|\alpha|}|\alpha|!\MM_{|\alpha|}
	\comma
	\alpha\in\mathbb{N}_0^{n+1}.
	\llabel{EQ112e}
\end{align}
	
	Recall that the radii in the doubling index of $Q\times I$ range over
	$0<r <\diam (Q\times I) = \sqrt{d} \ell (Q)$.
	The doubling inequality~\cite[Theorem~4.2]{DF1} provides constants
	$C,r_*>0$ such that 
	\begin{equation}
		\label{EQ41}
		\sup_{B_g(p,2r)}|\phi_\lambda|
		\leq
		e^{C(1+\sqrt\lambda)}
		\sup_{B_g(p,r)}|\phi_\lambda|,
	\end{equation}
	for any $p\in M$ and $0<r\leq r_*$.
	On compactly contained coordinate regions there exists a constant $L\geq 1$
	such that for
	$x\in Q$ and $p=\kappa^{-1}(x)$, we have
	\begin{align}
		\kappa\bigl(B_g(p,r/L)\bigr)
	\subset B_r(x)
		\quad
		\text{and}
		\quad
		\kappa^{-1}\bigl(B_{10dr}(x)\bigr)
		\subset B_g(p,10dLr).
		\llabel{EQ113a}
	\end{align}
Consequently, we arrive at
\begin{align}
	\frac{
		\sup_{B_{10dr} (x)}|\phi_\lambda|
	}{
		\sup_{B_r (x)}|\phi_\lambda|
	}
	\leq
	\frac{
		\sup_{B_g(p,10dLr)}|\phi_\lambda|
	}{
		\sup_{B_g(p,r/L)}|\phi_\lambda|
	}.
\llabel{EQ132a}
\end{align}

	We take an integer $J\geq1$ such that
	$2^J\geq10dL^2$.
	By $J$ iterations of~\eqref{EQ41}, we obtain
	\begin{align}
		\sup_{B_g(p,10dLr)}|\phi_\lambda|
		\leq
		e^{JC(1+\sqrt\lambda)}
		\sup_{B_g(p,r/L)}|\phi_\lambda|,
		\llabel{EQ133a}
	\end{align}
which leads to
	\begin{equation}
		\label{EQ46}
		\sup_{\substack{x\in Q\\0<r<\sqrt d\,\ell(Q)}}
		\log\frac{\sup_{B_{10dr} (x)}|\phi_\lambda|}
		{\sup_{B_r (x)}|\phi_\lambda|}
		\leq C(1+\sqrt\lambda).
	\end{equation}

	For $(x,t)\in Q\times I$, it is easy to see that
	\begin{align}
		\sup_{B_r (x,t)}|U_\lambda|
		\geq
		e^{\sqrt\lambda t}
		\sup_{B_r(x)}|\phi_\lambda|
		 \quad		
		 \text{and} 
		 \quad
		\sup_{B_{10dr} (x,t)}|U_\lambda|
		\leq
		e^{\sqrt\lambda(t+10dr)}
		\sup_{B_{10dr}(x)}|\phi_\lambda|.
		\llabel{EQ113d}
	\end{align}
	Combining these inequalities with
	\eqref{EQ46}, we obtain
	\begin{align}
		\log
		\frac{
			\sup_{B_{10dr} (x,t)}|U_\lambda|
		}{
			\sup_{B_r (x,t)}|U_\lambda|
		}
		\leq
		10dr\sqrt\lambda+C(1+\sqrt\lambda)
		\leq 
		C(1+\sqrt{\lambda}).
		\llabel{EQ113e}
	\end{align}
	Taking the supremum over all allowed
	centers and radii proves~\eqref{EQ44}.
	
	\subsection{The logarithmic weight}
	\label{App2}
	In this section, we verify that the logarithmic weight $\mathcal{M}_k= (\log (k+e))^{\beta k}$, where $0<\beta \leq 1$, satisfies the assumptions (i)--(v) in Section~\ref{sec:results} and the quasianalyticity condition \eqref{EQ105c}.
	
	It is clear that (i) holds.
	Let $L(x)=\log(x+e)$.
	The function $f(x)=\beta x\log L(x)$ is increasing on $[0,\infty)$
	and
	\begin{align}
	f''(x)=\frac{\beta((x+2e)L(x)-x)}{(x+e)^2L(x)^2}>0.
	\llabel{EQ111b}
	\end{align}
	Thus (ii)--(iii) hold.
	
	Let $j,k\in\mathbb{N}$.
	Concavity of $L$ and $L(0)=1$ imply
	\begin{align}
		\llabel{EQ120a}
		L(j)\geq
		\frac{j}{j+k} L(j+k)
		\quad
		\text{and}
		\quad
		L(k)
		\geq
		\frac{k}{j+k}
		L(j+k).
	\end{align} 
	Therefore, we arrive at
	\begin{align}
		\frac{\MM_{j+k}}{\MM_j\MM_k}
		\leq
		\left(\frac{j+k}{j}\right)^{j\beta}
		\left(\frac{j+k}{k}\right)^{k\beta}
		\leq2^{\beta(j+k)},	
		\llabel{EQ111a}
	\end{align}
	proving the assumption (v).  The assumption (iv) follows by taking
	$j=1$ in~\eqref{EQ105a}.
	
	Finally, we verify that \eqref{EQ105c} holds.
	For $k\in\mathbb{N}$, we have
	\begin{align}
	L(k)^\beta\leq\frac{\MM_k}{\MM_{k-1}}
	=L(k)^\beta\left(\frac{L(k)}{L(k-1)}\right)^{\beta(k-1)}.
	\label{EQ122a}
	\end{align}
	Using the Fundamental Theorem of Calculus and the concavity of $L$, we have
	\begin{align}
		L(x)-L(0)=
		\int_{0}^{x}
		L'(t)\,dt
		\geq 
		\int_{0}^{x}
		L'(x)\,dt
		=x L'(x)
		\comma
		x>0.
		\llabel{EQ121b}
	\end{align}
	Consequently, we obtain
	\begin{align}
		(\log L(x))'
		=
		\frac{L'(x)}{L(x)}
		\leq
		\frac{L(x)-1}{x L(x)}
		\leq
		\frac{1}{x}
		\comma x>0.
		\llabel{EQ121c}
	\end{align}
	For $k\geq 2$, another application of the Fundamental Theorem of Calculus gives
	\begin{align}
			\log L(k) - \log L(k-1)
			=
		\int_{k-1}^{k}
		(\log L(x))' \,dx
		\leq
		\int_{k-1}^{k}
		\frac{1}{x}\,dx
		=
		\log \left(
		1+\frac{1}{k-1}
		\right)
		\llabel{EQ121a}
	\end{align}
	which leads to
	\begin{align}
		\left(
		\frac{L(k)}{L(k-1)}
		\right)^{k-1}
		\leq
		\left(1+
		\frac{1}{k-1}
		\right)^{k-1}
		\leq e.
		\label{EQ122b}
	\end{align}
	Combining \eqref{EQ122a} with \eqref{EQ122b}, we get
	\begin{align}
		L(k)^\beta
		\leq 
	\frac{\mathcal{M}_k}{\mathcal{M}_{k-1}}	
		\leq
		e^\beta L(k)^\beta
		\comma
		k\in\mathbb{N}.
		\llabel{EQ122c}
	\end{align}
	Therefore, the series on the left-hand side of~\eqref{EQ105c} is divergent.
	
	\section*{Acknowledgments}
	\rm
	LL would like to thank Igor Kukavica for fruitful discussions.
	The author acknowledges the use of ChatGPT (GPT-6 Astra, OpenAI)
	in the preparation of this manuscript.
	All mathematical arguments
	and results were independently verified by the author, who takes full responsibility for the content of
	the paper.


\begin{thebibliography}{[LMNN]}

\bibitem[CZ]{CZ}
Y.~Cai and J.~Zhuge,
\emph{The sharp estimate of nodal sets for Dirichlet Laplace eigenfunctions in polytopes},
arXiv:2403.00279 (2024).

\bibitem[CKL]{CKL}
G.~Camliyurt, I.~Kukavica, and L.~Li,
\emph{Upper bounds of nodal sets for Gevrey regular parabolic equations},
Proc. Amer. Math. Soc.~\textbf{154} (2026), 4799--4809.

\bibitem[D]{D}
S.~Decio,
\emph{Hausdorff measure bounds for nodal sets of Steklov eigenfunctions},
Anal. PDE~\textbf{17} (2024), no.~4, 1237--1259.

\bibitem[Do]{Do}
R.~T.~Dong,
\emph{Nodal sets of eigenfunctions on Riemann surfaces},
J. Differential Geom.~\textbf{36} (1992), 493--506.

\bibitem[DF1]{DF1}
H.~Donnelly and C.~Fefferman,
\emph{Nodal sets of eigenfunctions on Riemannian manifolds},
Invent. Math.~\textbf{93} (1988), no.~1, 161--183.

\bibitem[DF2]{DF2}
H.~Donnelly and C.~Fefferman,
\emph{Nodal sets for eigenfunctions of the Laplacian on surfaces},
J. Amer. Math. Soc.~\textbf{3} (1990), no.~2, 333--353.

\bibitem[DF3]{DF3}
H.~Donnelly and C.~Fefferman,
\emph{Nodal sets of eigenfunctions: Riemannian manifolds with boundary},
in Analysis, et cetera, Academic Press, 1990, pp.~251--262.

\bibitem[Fu]{Fu}
S.~F\"urd\"os,
\emph{Geometric microlocal analysis in Denjoy--Carleman classes},
Pacific J. Math.~\textbf{307} (2020), no.~2, 303--351.

\bibitem[GR]{GR}
B.~Georgiev and G.~Roy-Fortin,
\emph{Polynomial upper bound on interior Steklov nodal sets},
J. Spectr. Theory~\textbf{9} (2019), no.~3, 897--919.

\bibitem[HS]{HS}
R.~Hardt and L.~Simon,
\emph{Nodal sets for solutions of elliptic equations},
J. Differential Geom.~\textbf{30} (1989), no.~2, 505--522.

\bibitem[He]{He}
H.~Hezari,
\emph{Upper bounds on the size of nodal sets for Gevrey and quasianalytic Riemannian manifolds},
Comm. Math. Phys.~\textbf{404} (2023), 1341--1359.

\bibitem[Hi]{Hi}
I.~I.~Hirschman, Jr.,
\emph{On the distributions of the zeros of functions belonging to certain quasi-analytic classes},
Amer. J. Math.~\textbf{72} (1950), no.~2, 396--406.

\bibitem[IK1]{IK1}
M.~Ignatova and I.~Kukavica,
\emph{Unique continuation and complexity of solutions to parabolic partial differential equations with Gevrey coefficients},
Adv. Differential Equations~\textbf{15} (2010), no.~9--10, 953--975.

\bibitem[IK2]{IK2}
M.~Ignatova and I.~Kukavica,
\emph{Strong unique continuation for higher order elliptic equations with Gevrey coefficients},
J. Differential Equations~\textbf{252} (2012), no.~4, 2983--3000.

\bibitem[KZZ]{KZZ}
C.~E.~Kenig, J.~Zhu, and J.~Zhuge,
\emph{Doubling inequalities and nodal sets in periodic elliptic homogenization},
Comm. Partial Differential Equations~\textbf{47} (2022), no.~3, 549--584.

\bibitem[K]{K}
I.~Kukavica,
\emph{Nodal volumes for eigenfunctions of analytic regular elliptic problems},
J. Anal. Math.~\textbf{67} (1995), 269--280.

\bibitem[KL]{KL}
I.~Kukavica and L.~Li,
\emph{Nodal set for the Schr\"odinger equation under a local growth condition},
Discrete Contin. Dyn. Syst. Ser. S~\textbf{26} (2026), 328--351.

\bibitem[L]{L}
F.~Lin,
\emph{Nodal sets of solutions of elliptic and parabolic equations},
Comm. Pure Appl. Math.~\textbf{44} (1991), no.~3, 287--308.

\bibitem[Lo1]{Lo1}
A.~Logunov,
\emph{Nodal sets of Laplace eigenfunctions: polynomial upper estimates of the Hausdorff measure},
Ann. of Math. (2)~\textbf{187} (2018), no.~1, 221--239.

\bibitem[Lo2]{Lo2}
A.~Logunov,
\emph{Nodal sets of Laplace eigenfunctions: proof of Nadirashvili's conjecture and of the lower bound in Yau's conjecture},
Ann. of Math. (2)~\textbf{187} (2018), no.~1, 241--262.

\bibitem[LM1]{LM1}
A.~Logunov and E.~Malinnikova,
\emph{Nodal sets of Laplace eigenfunctions: estimates of the Hausdorff measure in dimensions two and three},
in \emph{50 Years with Hardy Spaces}, Oper. Theory Adv. Appl., vol.~261,
Birkh\"auser/Springer, Cham, 2018, pp.~333--344.

\bibitem[LM2]{LM2}
A.~Logunov and E.~Malinnikova,
\emph{Review of Yau's conjecture on zero sets of Laplace eigenfunctions},
in \emph{Current Developments in Mathematics 2018},
Int. Press, Somerville, MA, 2020, pp.~179--212.

\bibitem[LMNN]{LMNN}
A.~Logunov, E.~Malinnikova, N.~Nadirashvili, and F.~Nazarov,
\emph{The sharp upper bound for the area of the nodal sets of Dirichlet Laplace eigenfunctions},
Geom. Funct. Anal.~\textbf{31} (2021), no.~5, 1219--1244.

\bibitem[LS]{LS}
F.~Lin and Z.~Shen,
\emph{Nodal sets and doubling conditions in elliptic homogenization},
Acta Math. Sin. (Engl. Ser.)~\textbf{35} (2019), no.~6, 815--831.

\bibitem[LTY]{LTY}
F.~Liu, L.~Tian, and X.~P.~Yang,
\emph{Measure upper bounds of nodal sets of Robin eigenfunctions},
Math. Z.~\textbf{306} (2024), no.~1, Paper No.~14, 14 pp.





\bibitem[LZ]{LZ}
F.~Lin and J.~Zhu,
\emph{Upper bounds of nodal sets for eigenfunctions of eigenvalue problems},
Math. Ann.~\textbf{382} (2022), no.~3--4, 1957--1984.

\bibitem[M]{M}
B.~S.~Mityagin,
\emph{The zero set of a real analytic function} (Russian),
Mat. Zametki~\textbf{107} (2020), no.~3, 473--475;
English transl., Math. Notes~\textbf{107} (2020), no.~3--4, 529--530.

\bibitem[NV]{NV}
A.~Naber and D.~Valtorta,
\emph{Volume estimates on the critical sets of solutions to elliptic PDEs},
Comm. Pure Appl. Math.~\textbf{70} (2017), no.~10, 1835--1897.

\bibitem[Y]{Y}
S.-T.~Yau (ed.),
\emph{Seminar on Differential Geometry},
Ann. of Math. Stud., vol.~102,
Princeton Univ. Press, Princeton, NJ, and University of Tokyo Press,
Tokyo, 1982.

\bibitem[Z1]{Z1}
J.~Zhu,
\emph{Doubling inequality and nodal sets for solutions of bi-Laplace equations},
Arch. Ration. Mech. Anal.~\textbf{232} (2019), no.~3, 1543--1595.

\bibitem[Z2]{Z2}
J.~Zhu, 
\emph{Upper bounds of nodal sets for solutions of bi-Laplace equations: II}, 
arXiv preprint arXiv:2603.04677, 2026.

\bibitem[ZZ]{ZZ}
J.~Zhu and J.~Zhuge,
\emph{Nodal sets of Dirichlet eigenfunctions in quasiconvex Lipschitz domains},
Adv. Math. 2026, ~\textbf{502}: 111150.

\end{thebibliography}
\end{document}